\documentclass[final,1p,times]{elsarticle}

\usepackage{amssymb}
 \usepackage{amsthm}
\usepackage{amscd}
\usepackage{amsmath}
\usepackage{amsfonts}
\usepackage{amssymb}
\usepackage{graphicx}

\numberwithin{equation}{section}

\newcommand{\ra}{\rightarrow}
\newcommand{\p}{\partial}
\newcommand{\f}{\frac}

\newcommand{\be}{\begin{equation}}
\renewcommand{\ra}{\rightarrow}
\newcommand{\ee}{\end{equation}}
\newcommand{\bea}{\begin{eqnarray}}
\newcommand{\eea}{\end{eqnarray}}
\newcommand{\bna}{\begin{eqnarray*}}
\newcommand{\ena}{\end{eqnarray*}}

\renewcommand{\le}{\left}
\newcommand{\ri}{\right}

\journal{***}

\begin{document}

\begin{frontmatter}

\title{Trudinger-Moser inequalities on the entire Heisenberg group}

\author{Yunyan Yang}
 \ead{yunyanyang@ruc.edu.cn}
\address{ Department of Mathematics,
Renmin University of China, Beijing 100872, P. R. China}

\begin{abstract}
Continuing our previous work (Cohn, Lam, Lu, Yang, Nonlinear
Analysis (2011), doi: 10.1016 /j.na.2011.09.053), we obtain a class
of Trudinger-Moser inequalities on the entire Heisenberg group,
which indicate what the best constants are. All the existing proofs
of similar inequalities on unbounded domain of the Euclidean space or the Heisenberg group are based on
 rearrangement argument. In this note, we propose a new approach to solve this problem.
Specifically we get the global Trudinger-Moser inequality by gluing local estimates with the help of
cut-off functions. Our method still works for similar problems when the Heisenberg group is replaced by the
Eclidean space or complete noncompact Riemannian manifolds.

\end{abstract}


\end{frontmatter}

%\tableofcontents
\section{Introduction}
Let $\mathbb{H}^n=\mathbb{R}^{2n}\times\mathbb{R}$ be the Heisenberg
group whose group action is defined by
\be\label{action}(x,y,t)\circ({x}^\prime,{y}^\prime,{t}^\prime)=(x+{x}^\prime,y+{y}^\prime,t+{t}^\prime+2(\langle
y,x^\prime\rangle- \langle x,{y}^\prime\rangle)),\ee where $x,y,
x^\prime, y^\prime \in\mathbb{R}^n$, $t,t^\prime\in\mathbb{R}$, and
$\langle\cdot,\cdot\rangle$ denotes the standard inner product in
$\mathbb{R}^n$. Let us denote the parabolic dilation in
$\mathbb{R}^{2n}\times\mathbb{R}$ by $\delta_\lambda$, namely,
$\delta_\lambda(\xi)=(\lambda x,\lambda y,\lambda^2t)$ for any
$\xi=(x,y,t)\in \mathbb{R}^{2n}\times\mathbb{R}$. The Jacobian
determinant of $\delta_\lambda$ is $\lambda^Q$, where $Q=2n+2$ is
the homogeneous dimension of $\mathbb{H}^n$.
 The following norm
 \be\label{norm}|\xi|_{h}=\le[\le(\sum_{i=1}^n(x_i^2+y_i^2)\ri)^2+t^2\ri]^{\f{1}{4}}\ee
 is homogeneous of degree one with respect to the dilation $\delta_\lambda$.
 The associated distance between two points $\xi$ and $\eta$ of $\mathbb{H}^n$ is defined accordingly by
 \be\label{distance}d_{h}(\xi,\eta)=|\eta^{-1}\circ \xi|_{h},\ee
 where $\eta^{-1}$ denotes the inverse of $\eta$ with respect to the group action, i.e. $\eta^{-1}=-\eta$.
 Obviously $d_h(\cdot,\cdot)$ is symmetric.
 The open ball of radius $r$ centered at $\xi$ is
 $$B_{h}(\xi,r)=\{\eta\in\mathbb{H}^n: d_{h}(\eta,\xi)<r\}.$$
 It is important to note that (see for example Stein \cite{Stein}, Section 5 of Chapter VIII)
 \be\label{volume}|B_{h}(\xi,r)|=|B_{h}(0,r)|=|B_{h}(0,1)|r^Q,\ee
 where $|\cdot|$ denotes the Lebesgue measure.
 The Lie algebra of
$\mathbb{H}^n$ is generated by the left-invariant vector fields
\be\label{1.5}T=\f{\p}{\p t},\,\,X_i=\f{\p}{\p x_i}+2y_i\f{\p}{\p t},\,\,Y_i=\f{\p}{\p y_i}-2x_i\f{\p}{\p t},\, i=1,\cdots,n.\ee
These generators satisfy the non-commutative formula
$[X_i,Y_i]=-4\delta_{ij}T$. Denote by $|\nabla_{\mathbb{H}^n}u|$ the
norm of the sub-elliptic gradient of a smooth function
$u:\mathbb{H}^n\ra\mathbb{R}$:
$$|\nabla_{\mathbb{H}^n}u|=\le(\sum_{i=1}^n\le((X_iu)^2+(Y_iu)^2\ri)\ri)^{{1}/{2}}.$$
Let $\Omega$ be an open set in $\mathbb{H}^n$.
We  use $W_0^{1,p}(\Omega)$ to denote the completion of
$C_0^\infty(\Omega)$ under the norm
\be\label{w-nor}\|u\|_{W_0^{1,p}(\Omega)}=\le(\int_{\Omega}\le(|\nabla_{\mathbb{H}^n}u|^p+|u|^p\ri)d\xi\ri)^{{1}/{p}}.\ee
In \cite{C-L}, Cohn-Lu proved a Trudinger-Moser inequality on
bounded smooth domains in the Hesenberg group $\mathbb{H}^n$.
Precisely, there exists some constant $C_n$ depending only on $n$
such that for all bounded smooth domain $\Omega\subset\mathbb{H}^n$,
if $u\in W_0^{1,Q}(\Omega)$ satisfies
$\|\nabla_{\mathbb{H}^n}u\|_{L^Q(\Omega)}\leq 1$, then
\be\label{CL}\int_\Omega
 e^{\alpha_Q|u|^{Q^\prime}}d\xi\leq C_n|\Omega|,\ee
 where $Q^\prime=Q/(Q-1)$, $\alpha_Q=Q\sigma_{Q}^{1/(Q-1)}$,
 $\sigma_Q=\Gamma(\frac{1}{2})\Gamma(n+\frac{1}{2})\omega_{2n-1}/n!$,
 $\omega_{2n-1}$ is the surface area of the unit sphere in
 $\mathbb{R}^{2n}$. Furthermore, the
 integrals of all  $u\in W_0^{1,Q}(\Omega)$ satisfying
$\|\nabla_{\mathbb{H}^n}u\|_{L^Q(\Omega)}\leq 1$ are not uniformly
bounded if $\alpha_Q$ is replaced by any larger number. Recently,
Cohn, Lam, Lu and the author \cite{CLLY} obtained a Trudinger-Moser
inequality on the Heisenberg group $\mathbb{H}^n$. Note that
$W^{1,Q}(\mathbb{H}^n)$ is the completion of
$C_0^\infty(\mathbb{H}^n)$ under the norm (\ref{w-nor}) with
$\Omega$ replaced by $\mathbb{H}^n$.
We have the following:\\

\noindent{\bf Theorem A} (\cite{CLLY}). {\it There exists some constant $\alpha^\ast: 0<\alpha^\ast\leq \alpha_Q$ such that
 for any pair $\beta$ and $\alpha$ satisfying $0\leq\beta<Q$,
$0<\alpha\leq \alpha^\ast$, and
$\f{\alpha}{\alpha^\ast}+\f{\beta}{Q}\leq 1$, there holds
\be\label{TMH}\sup_{\|u\|_{W^{1,Q}(\mathbb{H}^n)}\leq
1}\int_{\mathbb{H}^n}
\f{1}{|\xi|_h^\beta}\le\{e^{\alpha|u|^{Q\,^\prime}}-\sum_{k=0}^{Q-2}
\f{\alpha^k|u|^{kQ\,^\prime}}{k!}\ri\}d\xi<\infty.\ee
 When $\f{\alpha}{\alpha^\ast}+\f{\beta}{Q}> 1$, the integral in $(\ref{TMH})$ is still finite
 for any $u\in W^{1,Q}(\mathbb{H}^n)$, but the supremum is
 infinite if further $\f{\alpha}{\alpha_Q}+\f{\beta}{Q}>1$.}\\

 Theorem A is an analogue of (Adimurthi-Yang \cite{Adi-Yang}, Theorem 1.1).
  Earlier works on this topic (Trudinger-Moser inequalities on unbounded domain of $\mathbb{R}^n$) were done by Cao \cite{Cao}, Panda \cite{Panda}, do \'O \cite{doo},
 Ruf \cite{Ruf}, Li-Ruf \cite{Li-Ruf} and others. The proof of Theorem A is based on
 symmetrization argument, radial lemma and the Young inequality.
 Note that $\alpha^*$ in Theorem A is not explicitly known. A natural question is what the best constant $\alpha$ for (\ref{TMH}) is.
 Denote an equivalent norm in $W^{1,Q}(\mathbb{H}^n)$ by
 \be\label{1tau}\|u\|_{1,\tau}=\le(\int_{\mathbb{H}^n}(|\nabla_{\mathbb{H}^n}u|^Q+\tau|u|^Q)d\xi\ri)^{\f{1}{Q}}\ee
 for any fixed number $\tau>0$. Our main result is the following:\\

 \noindent{\bf Theorem 1.1.} {\it Let $\tau$ be any positive real number. Let $Q$, $Q^\prime$
 and $\alpha_Q$ be as in (\ref{CL}).
 For any $\beta:0\leq\beta<Q$ and $\alpha:0<\alpha<\alpha_Q(1-\beta/Q)$,
 there holds
 \be\label{Ttau}\sup_{\|u\|_{1,\tau}\leq
1}\int_{\mathbb{H}^n}
\f{1}{|\xi|_h^\beta}\le\{e^{\alpha|u|^{Q\,^\prime}}-\sum_{k=0}^{Q-2}
\f{\alpha^k|u|^{kQ\,^\prime}}{k!}\ri\}d\xi<\infty.\ee
When $\alpha>\alpha_Q(1-\beta/Q)$, the above integral is still finite
 for any $u\in W^{1,Q}(\mathbb{H}^n)$, but the supremum is
 infinite.} \\

 Clearly Theorem 1.1 implies that the best constant for the inequality (\ref{Ttau}) is $\alpha_Q(1-\beta/Q)$. But we do not know
 whether or not (\ref{Ttau}) still holds when $\alpha=\alpha_Q(1-\beta/Q)$. Even so, $(\ref{Ttau})$ gives more
 information than (\ref{TMH}).

  According to the author's knowledge, the existing proofs of Trudinger-Moser inequalities for unbounded domains
  are all based on the rearrangement theory \cite{H-L}.
  It is not known that whether or not this technique can be successfully applied to the Heisenberg group case.

  To prove Theorem 1.1,
    we propose a new approach. The idea can be described as follows. Firstly, using (\ref{CL}),
  we derive a local Trudinger-Moser inequality, namely, for any fixed $r>0$
  and all $\xi_0\in\mathbb{H}^n$, there exists some constant $C$ depending only on $n$, $r$ and $\beta$ such that
  \be\label{loc1}\int_{B_h(\xi_0,r)}\f{1}{|\xi|_h^\beta}\le\{e^{\alpha|u|^{Q\,^\prime}}-\sum_{k=0}^{Q-2}
\f{\alpha^k|u|^{kQ\,^\prime}}{k!}\ri\}d\xi\leq C\int_{B_h(\xi_0,r)}|\nabla_{\mathbb{H}^n}u|^Qd\xi\ee
provided that $0\leq \alpha<\alpha_Q(1-\beta/Q)$ and $\int_{B_h(\xi_0,r)}|\nabla_{\mathbb{H}^n}u|^Qd\xi\leq 1$. Secondly, fixing sufficiently large $r>0$, we
 select a specific sequence of Heisenberg balls $\{B_h(\xi_i,r)\}_{i=1}^\infty$ to cover the Heisenberg group $\mathbb{H}^n$. Then
 we choose appropriate cut-off function $\phi_i$ on each $B_h(\xi_i,r)$. Finally, we obtain (\ref{Ttau}) by gluing all local estimates (\ref{loc1}) for
 $\phi_iu$. We remark that our method still works for similar problems when the Heisenberg group is replaced by the
Eclidean space or complete noncompact Riemannian manifolds. In the Eclidean space case, $\tau$ can also be arbitrary in (\ref{Ttau}). But in the manifold
case, the choice of $\tau$ may depend on the geometric structure (see \cite{Yang-manifolds}, Theorem 2.3).
 As an easy consequence of Theorem 1.1 (in fact a special case $\beta=0$), the following corollary holds.\\

 \noindent{\bf Corollary 1.2.} {\it Let $Q=2n+2$. For any $q\geq Q$, $W^{1,Q}(\mathbb{H}^n)$ is continuously embedded in
 $L^q(\mathbb{H}^n)$.}\\

   The remaining part of this note is organized as follows. In
   section 2, we prove a covering lemma for $\mathbb{H}^n$;
   Cut-off functions are selected for the subsequent analysis in
   section 3; The proof of Theorem 1.1 is completed in section 4.

\section{A covering lemma for the Heisenberg group}

In this section, we will use a sequence of Heisenberg balls with the same radius to cover the entire Heisenberg group $\mathbb{H}^n$. We require these balls
to satisfy the following properties: $(i)$ For any $\xi\in \mathbb{H}^n$, $\xi$ belongs to at most $N$ balls for some constant integer $N$ which is
independent of the base point $\xi$; $(ii)$ If the radius of those balls becomes appropriately smaller, then they are disjoint.

Firstly, we need to understand the Heisenberg distance between two
points of the Heisenberg group $\mathbb{H}^n$. The following
two properties are more or less standard. We prefer to present them by our own way.\\

\noindent{\bf Proposition 2.1.} {\it Let $\xi$ and $\eta$ be two points of $\mathbb{H}^n$. There holds
$$|\eta^{-1}\circ \xi|_h\leq 3(|\xi|_h+|\eta|_h),$$
where $|\cdot|_h$ is the homogeneous norm defined by (\ref{norm}).}\\

\noindent{\it Proof.} Write $\xi=(x,y,t)$, $\eta=(x^\prime,y^\prime,t^\prime)$. Then (\ref{action}) gives
$$\eta^{-1}\circ\xi=(x-x^\prime,y-y^\prime,t-t^\prime-2(\langle y,x^\prime\rangle-\langle x,y^\prime\rangle)).$$
Since $(|x-x^\prime|^2+|y-y^\prime|^2)^{1/2}\leq (|x|^2+|y|^2)^{1/2}+(|x^\prime|^2+|y^\prime|^2)^{1/2}$ and
$$\le|2\le(\langle y,x^\prime\rangle-\langle x,y^\prime\rangle\ri)\ri|\leq |x|^2+|y|^2+|x^\prime|^2
+|y^\prime|^2,$$
we have by using the inequality
$\sqrt{a+b}\leq \sqrt{a}+\sqrt{b}$ ($a\geq 0$, $b\geq 0$) repeatedly
\bna
 |\eta^{-1}\circ\xi|_h&=&\le[\le(\sum_{i=1}^n\le((x_i-x_i^\prime)^2+(y_i-y_i^\prime)^2\ri)\ri)^2+(t-t^\prime
 -2\le(\langle y,x^\prime\rangle-\langle x,y^\prime\rangle)\ri)^2\ri]^{\f{1}{4}}\\
  &\leq&\le(\sum_{i=1}^n\le((x_i-x_i^\prime)^2+(y_i-y_i^\prime)^2\ri)\ri)^{\f{1}{2}}
 +\le|t-t^\prime
 -2\le(\langle y,x^\prime\rangle-\langle x,y^\prime\rangle\ri)\ri|^{\f{1}{2}}\\
 &\leq&2\le(\sum_{i=1}^n\le(x_i^2+y_i^2\ri)\ri)^{\f{1}{2}}+2\le(\sum_{i=1}^n\le({x_i^\prime}^2+{y_i^\prime}^2\ri)\ri)^{\f{1}{2}}
 +|t|^{\f{1}{2}}+|t^\prime|^{\f{1}{2}}\\
 &\leq& 3(|\xi|_h+|\eta|_h).
\ena
$\hfill\Box$

\noindent{\bf Proposition 2.2.} {\it Let $\xi$, $\eta$, $\zeta$ be arbitrary points of $\mathbb{H}^n$. Then we have
$$d_h(\xi,\eta)\leq 3\le(d_h(\xi,\zeta)+d_h(\zeta,\eta)\ri),$$
where $d_h(\cdot,\cdot)$ is the distance function defined by (\ref{distance}).}\\

\noindent{\it Proof.} Note that $|\gamma^{-1}|_h=|\gamma|_h$ for all $\gamma\in\mathbb{H}^n$. It follows from Proposition 2.1 that
\bna
 d_h(\xi,\eta)&=&|\eta^{-1}\circ\xi|_h\\
 &=&|\eta^{-1}\circ\zeta\circ\zeta^{-1}\circ\xi|_h\\
 &\leq&3(|\eta^{-1}\circ\zeta|_h+|\zeta^{-1}\circ\xi|_h)\\
 &=&3\le(d_h(\zeta,\eta)+d_h(\xi,\zeta)\ri).
\ena
This gives the desired result. $\hfill\Box$\\

Secondly, by adapting an argument of (Hebey \cite{Hebey}, Lemma 1.6), we obtain the following useful covering lemma.\\

\noindent{\bf Lemma 2.3.} {\it Let $\rho>0$ be given. There exists a sequence $(\xi_i)$ of points of $\mathbb{H}^n$
such that for any $r\geq \rho$:\\
$(i)$ $\cup_{i}B_h(\xi_i,\rho)=\mathbb{H}^n$ and for any $i\not=j$,
$B_h(\xi_i,\rho/6)\cap
B_h(\xi_j,\rho/6)=\varnothing$;\\
$(ii)$ for any $\xi\in\mathbb{H}^n$, $\xi$ belongs to at most
$[(24r/\rho)^Q]$ balls $B_h(\xi_i,r)$,
where $[(24r/\rho)^Q]$ denotes the integral part of $(24r/\rho)^Q$.}\\

\noindent{\it Proof.} Firstly, we {\it claim} that there exists a
sequence $(\xi_i)$ of points of $\mathbb{H}^n$ such that
\be\label{cov}\cup_{i}B_h(\xi_i,\rho)=\mathbb{H}^n\,\, {\rm and}\,\,
\forall i\not=j, B_h(\xi_i,\rho/6)\cap
B_h(\xi_j,\rho/6)=\varnothing.\ee To see this, we set
$$X_\rho=\le\{{\rm sequence}\,\, (\xi_i)_{i\in I}:\xi_i\in \mathbb{H}^n, I\,\,{\rm is\,\,countable\,\,and}\,\,\forall i\not=j,
d_h(\xi_i,\xi_j)\geq \rho\ri\}.$$
Then $X_\rho$ is partially ordered by inclusion and every element in $X_\rho$ has an upper bound in the sense of inclusion.
Hence, by Zorn's lemma, $X_\rho$ contains a maximal element $(\xi_i)_{i\in I}$. On one hand, if
$\cup_{i}B_h(\xi_i,\rho)\not=\mathbb{H}^n$, then there exists a point $\xi\in\mathbb{H}^n$ such that $d_h(\xi_i,\xi)\geq \rho$ for
all $i\in I$. This contradicts the maximality of $(\xi_i)_{i\in I}$. Hence  $\cup_{i}B_h(\xi_i,\rho)=\mathbb{H}^n$. On the other hand,
if $B_h(\xi_i,\rho/6)\cap B_h(\xi_j,\rho/6)\not=\varnothing$ for some $i\not=j$, then we can take some
$\eta\in B_h(\xi_i,\rho/6)\cap B_h(\xi_j,\rho/6)$. It follows from Proposition 2.2 that
\bna
d_h(\xi_i,\xi_j)&\leq& 3\le(d_h(\xi_i,\eta)+d_h(\eta,\xi_j)\ri)\\
&<&3\le(\f{\rho}{6}+\f{\rho}{6}\ri)=\rho.
\ena
This contradicts the fact that $d_h(\xi_i,\xi_j)\geq \rho$ for any $i\not=j$. Thus our claim (\ref{cov}) holds.

Assume $(\xi_i)$ satisfies (\ref{cov}). For any fixed $r>0$ and $\xi\in\mathbb{H}^n$ we set
$$I_r(\xi)=\le\{i\in I: \xi\in B_h(\xi_i,r)\ri\}.$$
By (\ref{volume}) and Proposition 2.2, we have for $r\geq \rho$
\bna
|B_h(\xi,r)|&=&4^{-Q}|B_h(\xi,4r)|\\
&\geq&4^{-Q}\sum_{i\in I_r(\xi)}|B_h(\xi_i,\rho/6)|\\
&=&4^{-Q}\,\,{\rm Card}\,\,I_r(\xi)\,\,({\rho}/{6})^Q|B_h(0,1)|,
\ena
where ${\rm Card}\,\,I_r(\xi)$ denotes the cardinality of the set $I_r(\xi)$. As a consequence, for $r\geq \rho$ there holds
$${\rm Card}\,\,I_r(\xi)\leq (24r/\rho)^Q.$$
This completes the proof of the lemma. $\hfill\Box$

\section{Cut-off functions on Heisenberg balls}

In this section, we will construct cut-off functions on Heisenberg balls. To do this, we first estimate the gradient of the distance function as follows.\\

\noindent {\bf Lemma 3.1.} {\it Let $\xi_0$ be any fixed point of $\mathbb{H}^n$. Define a function $\rho(\xi)=d_h(\xi,\xi_0)$. Then we have
$|\nabla_{\mathbb{H}^n}\rho(\xi)|\leq 1$ for any $\xi\not=\xi_0$.}\\

\noindent{\it Proof.} Write $\xi=(x_1,\cdots,x_n,y_1,\cdots,y_n,t)$ and $\xi_0=(x_{01},\cdots,x_{0n},y_{01},\cdots,y_{0n},t_0)$.
For any $\xi\not=\xi_0$, we set
$$E=\sum_{i=1}^n\le((x_i-x_{0i})^2+(y_i-y_{0i})^2\ri),\quad F=t-t_0-2\sum_{i=1}^n(x_iy_{0i}-y_ix_{0i}).$$
Then by (\ref{action}) and (\ref{distance}),
$$\rho(\xi)=|\xi_0^{-1}\circ\xi|_h=\le(E^2+F^2\ri)^{1/4}.$$
We calculate
\bna
\f{\p}{\p x_i}\rho=\rho^{-3}\le((x_i-x_{0i})E-y_{0i}F\ri),\quad
2y_i\f{\p}{\p t}\rho=\rho^{-3}y_iF,
\ena
and then by (\ref{1.5}),
$$X_i\rho=\f{\p}{\p x_i}\rho+2y_i\f{\p}{\p t}\rho=\rho^{-3}\le((x_i-x_{0i})E+(y_i-y_{0i})F\ri).$$
Similarly we have
$$\f{\p}{\p y_i}\rho=\rho^{-3}\le((y_i-y_{0i})E+x_{0i}F\ri)$$
and thus by (\ref{1.5}),
$$Y_i\rho=\f{\p}{\p y_i}\rho-2x_i\f{\p}{\p t}\rho=\rho^{-3}\le((y_i-y_{0i})E+(x_{0i}-x_i)F\ri).$$
It follows that
$$(X_i\rho)^2+(Y_i\rho)^2=\rho^{-6}\le((y_i-y_{0i})^2+(x_i-x_{0i})^2\ri)(E^2+F^2).$$
Note that $E^2+F^2=\rho^4$. We obtain
\bna
|\nabla_{\mathbb{H}^n}\rho|&=&\le(\sum_{i=1}^n\le((X_i\rho)^2+(Y_i\rho)^2\ri)\ri)^{1/2}\\
&=&\rho^{-3}E^{1/2}(E^2+F^2)^{1/2}\\
&=&\rho^{-1}E^{1/2}\leq 1.
\ena
This completes the proof of the lemma. $\hfill\Box$\\

Now we construct cut-off functions. Let $\phi:\mathbb{R}\ra\mathbb{R}$ be a smooth function such that $0\leq \phi\leq 1$, $\phi\equiv 1$ on
the interval $[-1,1]$, $\phi\equiv 0$ on $(-\infty,-2)\cup(2,\infty)$, and $|\phi^\prime(t)|\leq 2$ for all $t\in \mathbb{R}$.
Let $r>0$ be given.
Define a function on $\mathbb{H}^n$ by
\be\label{cut}\phi_{0}(\xi)=\phi\le(\f{d_h(\xi,\xi_0)}{r}\ri).\ee
Then $\phi_0$ is a cut-off function supported on the Heisenberg ball $B_h(\xi_0,2r)$. The estimate of the gradient of $\phi_0$ is very
important for the subsequent analysis. Precisely we have the following:\\

\noindent{\bf Lemma 3.2.} {\it For any fixed $r>0$ and $\xi_0\in\mathbb{H}^n$, let $\phi_0$ be defined by (\ref{cut}).
Then $\phi_0$ is supported in $B_{h}(\xi_0,2r)$, $0\leq\phi_0\leq 1$,
$\phi_0\equiv 1$ on $B_h(\xi_0,r)$, and
$|\nabla_{\mathbb{H}^n}\phi_0(\xi)|\leq 2/r$ for all $\xi\in\mathbb{H}^n$.}\\

\noindent{\it Proof.} We only need to explain the last assertion, namely
$|\nabla_{\mathbb{H}^n}\phi_0(\xi)|\leq 2/r$ for all $\xi\in\mathbb{H}^n$. Since $\phi_0\equiv 1$ on $B_h(\xi_0,r)$,
we have $\nabla_{\mathbb{H}^n}\phi_0\equiv 0$ on $B_h(\xi_0,r)$, particularly $\nabla_{\mathbb{H}^n}\phi_0(0)= 0$.
For $\xi\not= \xi_0$, a simple calculation shows
$$\nabla_{\mathbb{H}^n}\phi_0(\xi)=\f{1}{r}\phi^\prime\nabla_{\mathbb{H}^n}d_h(\xi,\xi_0).$$
This together with Lemma 3.1 and $|\phi^\prime|\leq 2$ concludes the last assertion. $\hfill\Box$

\section{Proof of Theorem 1.1}

In this section, we will prove Theorem 1.1. For simplicity, we define a
smooth function $\zeta:\mathbb{N}\times\mathbb{R}\ra \mathbb{R}$ by
\be\label{MF}\zeta(m,s)=e^s-\sum_{k=0}^{m-2}\f{s^k}{k!},\quad\forall m\geq 2.\ee
As we promised in the introduction, we first derive a local Trudinger-Moser inequality  for
the Heisenberg group $\mathbb{H}^n$ by using (\ref{CL}). Let $Q$, $Q^\prime$ and $\alpha_Q$ be
given by (\ref{CL}). Then we have the following:\\

\noindent{\bf Lemma 4.1.}
{\it Let $r>0$ be given and $\xi_0$ be any point of $\mathbb{H}^n$. If $0\leq \beta<Q$, $0\leq \alpha\leq\alpha_Q(1-\beta/Q)$,
and $w\in W_0^{1,Q}(B_h(\xi_0,r))$ satisfies $\int_{B_h(\xi_0,r)}|\nabla_{\mathbb{H}^n}w|^Qd\xi\leq 1$, then there exists
 some constant $C$ depending only on $n$, $r$ and $\beta$ such that
  \be\label{loc}\int_{B_h(\xi_0,r)}\f{1}{|\xi|_h^\beta}\zeta(Q,\alpha|w|^{Q^\prime})d\xi\leq C\int_{B_h(\xi_0,r)}|\nabla_{\mathbb{H}^n}w|^Qd\xi.\ee}

\noindent{\it Proof.} Using Proposition 2.2, we have that
$$|\xi_0|_h\leq 3(d_h(\xi,\xi_0)+|\xi|_h),\quad\forall \xi\in\mathbb{H}.$$
If $|\xi_0|_h>6r$, then for any $\xi\in B_h(\xi_0,r)$ there holds
\be\label{gq}|\xi|_h\geq \f{|\xi_0|_h}{3}-d_h(\xi,\xi_0)>r.\ee
Let $\widetilde{w}=w/\|\nabla_{\mathbb{H}^n}
 w\|_{L^Q({B}_h(\xi_0,r))}$. Since $\|\nabla_{\mathbb{H}^n}
 w\|_{L^Q({B}_h(\xi_0,r))}\leq 1$ and $0\leq\alpha\leq\alpha_Q(1-\beta/Q)$, we have
 \bea
\zeta\le(Q,\alpha|w|^{Q^\prime}\ri)&=&\sum_{k=Q-1}^\infty
 \f{\alpha^k|w|^{Q^\prime k}}{k!}{\nonumber}\\
 &=&\sum_{k=Q-1}^\infty
 \f{\alpha^k\|\nabla_{\mathbb{H}^n}
 w\|_{L^Q({B}_h(\xi_0,r))}^{Q^\prime k}|\widetilde{w}|^{Q^\prime k}}{k!}{\nonumber}\\
 \label{4}
 &\leq&\|\nabla_{\mathbb{H}^n}
 w\|_{L^Q({B}_h(\xi_0,r))}^Q\zeta\le(Q,\alpha|\widetilde{w}|^{Q^\prime}\ri).
 \eea
 By (\ref{volume}) and (\ref{CL}),
 $$\int_{B_h(\xi_0,r)}\zeta\le(Q,\alpha_Q|\widetilde{w}|^{Q^\prime}\ri)d\xi\leq
 C_nr^Q|B_h(0,1)|,$$
 where $C_n$ is given by (\ref{CL}). Hence when $|\xi_0|_h>6r$ and $0\leq\alpha\leq \alpha_Q(1-\beta/Q)$,
 we have by using (\ref{gq}) and (\ref{4}),
 \bna
 \int_{B_h(\xi_0,r)}\f{1}{|\xi|_h^\beta}\zeta\le(Q,\alpha|w|^{Q^\prime}\ri)d\xi&\leq&r^{-\beta}
 \int_{B_h(\xi_0,r)}\zeta\le(Q,\alpha|w|^{Q^\prime}\ri)d\xi\\
 &\leq&C_nr^{Q-\beta}|B_h(0,1)|\int_{B_h(\xi_0,r)}|\nabla_{\mathbb{H}^n}w|^Qd\xi.
 \ena

 In the following we assume $|\xi_0|_h\leq 6r$. If $\xi\in B_h(\xi_0,r)$, then Proposition 2.2 implies that
 $$|\xi|_h\leq 3(d_h(\xi,\xi_0)+|\xi_0|_h)< 21r.$$
 H\"older's inequality together with (\ref{CL}) implies that there exits some constant $\widetilde{C}$ depending only on
 $n$, $r$ and $\beta$ such that
 $$\int_{B_h(\xi_0,r)}\f{1}{|\xi|_h^\beta}\zeta\le(Q,\alpha|\widetilde{w}|^{Q^\prime}\ri)d\xi
 \leq\int_{|\xi|_h\leq 21r}\f{1}{|\xi|_h^\beta}\zeta\le(Q,\alpha|\widetilde{w}|^{Q^\prime}\ri)d\xi\leq \widetilde{C}.$$
 It then follows from (\ref{4}) that
 \bna
 \int_{B_h(\xi_0,r)}\f{1}{|\xi|_h^\beta}\zeta\le(Q,\alpha|w|^{Q^\prime}\ri)d\xi&\leq&
 \|\nabla_{\mathbb{H}^n}
 w\|_{L^Q({B}_h(\xi_0,r))}^Q\int_{B_h(\xi_0,r)}\f{1}{|\xi|_h^\beta}\zeta\le(Q,\alpha|\widetilde{w}|^{Q^\prime}\ri)d\xi\\
 &\leq&\widetilde{C}\int_{B_h(\xi_0,r)}|\nabla_{\mathbb{H}^n}w|^Qd\xi.
 \ena
   Hence (\ref{loc}) holds. $\hfill\Box$\\

 {\it Proof of Theorem 1.1.} Firstly, we prove (\ref{Ttau}). Let $\tau>0$ and $\alpha: 0\leq \alpha<\alpha_Q(1-\beta/Q)$ be fixed.
 Since $C_0^\infty(\mathbb{H}^n)$ is dense in $W^{1,Q}(\mathbb{H}^n)$ under the norm (\ref{1tau}), 
 it suffices to prove (\ref{Ttau}) for all $u\in C_0^\infty(\mathbb{H}^n)$ with
 \be\label{nm}\int_{\mathbb{H}^n}(|\nabla_{\mathbb{H}^n}u|^Q+\tau|u|^Q)d\xi\leq 1.\ee
  Assume $u\in C^\infty(\mathbb{H}^n)$ satisfies (\ref{nm}). Let $r>0$ be a sufficiently large number to be determined later.
 By Lemma 2.3, there exists a sequence $(\xi_i)$ of points of $\mathbb{H}^n$ such that
 \be\label{cv}\cup_iB_h(\xi_i,r)=\mathbb{H}^n\,\,{\rm and}\,\,\forall
 i\not=j,\,\,
 B_h(\xi_i,r/6)\cap B_h(\xi_j,r/6)=\varnothing,\ee
 and for any $\xi\in\mathbb{H}^n$,
 \be\label{4.7}\xi\,\,{\rm belongs\,\, to\,\, at\,\, most}\,\, 48^Q \,\,{\rm balls}\,\, B_h(\xi_i,2r).\ee

 Let $\phi$ be a smooth function given by (\ref{cut}). For each $\xi_i$, we set
 $$\phi_i(\xi)=\phi\le(\f{d_h(\xi,\xi_i)}{r}\ri),\quad\forall \xi\in\mathbb{H}^n.$$
  It follows from Lemma 3.2 that
  $0\leq \phi_i\leq 1$, $\phi_i\equiv 1$ on $B_h(\xi_i,r)$,
  $\phi_i\equiv 0$ outside $B_h(\xi_i,2r)$, and
  \be\label{gd}|\nabla_{\mathbb{H}^n}\phi_i(\xi)|\leq \f{2}{r},
  \,\,\forall \xi\in\mathbb{H}^n.  \ee
  Clearly $\phi_i^2 u\in W_0^{1,Q}\le(B_h(\xi_i,2r)\ri)$. Since $u$ satisfies (\ref{nm}), we have that
  $$\int_{\mathbb{H}^n}|\nabla_{\mathbb{H}^n}u|^Qd\xi\leq 1,\,\,{\rm and}\,\,
  \int_{\mathbb{H}^n}|u|^Qd\xi\leq \f{1}{\tau}.$$
  Minkowski inequality together with (\ref{gd}) and $0\leq \phi_i\leq 1$ leads to
 \bea
 \le(\int_{B_h(\xi_i,2r)}|\nabla_{\mathbb{H}^n}(\phi_i^2u)|^Qd\xi\ri)^{1/Q}&\leq&
 \le(\int_{B_h(\xi_i,2r)}\phi_i^{2Q}|\nabla_{\mathbb{H}^n}u|^Qd\xi\ri)^{1/Q}+
 \le(\int_{B_h(\xi_i,2r)}|\nabla_{\mathbb{H}^n}\phi_i^2|^Q|u|^Qd\xi\ri)^{1/Q}{\nonumber}\\
 &\leq&\le(\int_{B_h(\xi_i,2r)}|\nabla_{\mathbb{H}^n}u|^Qd\xi\ri)^{1/Q}+
 \f{4}{r}\le(\int_{B_h(\xi_i,2r)}|u|^Qd\xi\ri)^{1/Q}{\nonumber}\\\label{M-K}
 &\leq&1+\f{4}{\tau r}.
 \eea
 Define
 $\widetilde{u}_i={\phi_i^2u}/(1+\f{4}{\tau r})$.
 Then $\widetilde{u}_i\in W_0^{1,Q}(B_h(\xi_i,2r))$ and $\int_{B_h(\xi_i,2r)}
 |\nabla_{\mathbb{H}^n}\widetilde{u}_i|^Qd\xi\leq 1$. Since $\alpha<\alpha_Q(1-\beta/Q)$,
 we can select $r$ sufficiently large such that
 $$\alpha\le(1+\f{4}{\tau r}\ri)^{Q^\prime}<\alpha_Q(1-\beta/Q).$$
 This together with Lemma 4.1 implies that there exists some
 constant $C$ depending only on $n$, $r$ and $\beta$ such that
 \bea\int_{B_h(\xi_i,2r)}\f{1}{|\xi|_h^\beta}\zeta
 \le(Q,\alpha|\phi_i^2u|^{Q^\prime}\ri)d\xi&=&\int_{B_h(\xi_i,2r)}\f{1}{|\xi|_h^\beta}\zeta
 \le(Q,\alpha\le(1+\f{4}{\tau r}\ri)^{Q^\prime}|\widetilde{u}_i|^{Q^\prime}\ri)d\xi{\nonumber}\\
 &\leq& C\int_{B_h(\xi_i,2r)}|\nabla_{\mathbb{H}^n}\widetilde{u}_i|^Qd\xi.\nonumber\\
 &\leq&C\int_{B_h(\xi_i,2r)}|\nabla_{\mathbb{H}^n}(\phi_i^2u)|^Qd\xi.\label{l-e}\eea
 Combining (\ref{cv}) and (\ref{l-e}), we obtain
 \bea
  \int_{\mathbb{H}^n}\f{1}{|\xi|_h^\beta}\zeta\le(Q,\alpha|u|^{Q^\prime}\ri)d\xi&\leq&\sum_{i}
 \int_{B_h(\xi_i,r)}\f{1}{|\xi|_h^\beta}\zeta\le(Q,\alpha|\phi_i^2u|^{Q^\prime}\ri)d\xi{\nonumber}\\{\nonumber}
 &\leq&\sum_{i}
 \int_{B_h(\xi_i,2r)}\f{1}{|\xi|_h^\beta}\zeta\le(Q,\alpha|\phi_i^2u|^{Q^\prime}\ri)d\xi\\\label{16}
 &\leq&C\sum_{i}\int_{\mathbb{H}^n}|\nabla(\phi_i^2 u)|^Qd\xi.
 \eea
 Using the inequality $|a+b|^Q\leq 2^Q|a|^Q+2^Q|b|^Q$, $\forall a,b\in\mathbb{R}$,
 $0\leq \phi_i\leq 1$ and  (\ref{gd}), we get
  \bna
  \int_{\mathbb{H}^n}|\nabla_{\mathbb{H}^n}(\phi_i^2 u)|^Qd\xi&\leq&
  2^Q\int_{\mathbb{H}^n}\le(\phi_i^{2Q}|\nabla_{\mathbb{H}^n}
  u|^Q+|\nabla_{\mathbb{H}^n}\phi_i^2|^Q|u|^Q\ri)d\xi{\nonumber}\\\label{26}
  &\leq&2^Q\int_{\mathbb{H}^n}\phi_i|\nabla_{\mathbb{H}^n}
  u|^Qd\xi+\le(\f{8}{r}\ri)^Q\int_{\mathbb{H}^n}\phi_i|u|^Qd\xi.
 \ena
 In view of (\ref{4.7}), it then follows that
 \bna
  \sum_i\int_{\mathbb{H}^n}|\nabla_{\mathbb{H}^n}(\phi_i^2 u)|^Qd\xi&\leq&
  2^Q\sum_i\int_{\mathbb{H}^n}\phi_i|\nabla_{\mathbb{H}^n}
  u|^Qd\xi+\le(\f{8}{r}\ri)^Q\sum_i\int_{\mathbb{H}^n}\phi_i|u|^Qd\xi
  \\&\leq&96^Q\int_{\mathbb{H}^n}|\nabla_{\mathbb{H}^n}
  u|^Qd\xi+\le(\f{384}{r}\ri)^Q\int_{\mathbb{H}^n}|u|^Qd\xi.
   \ena
 This together with (\ref{16}) implies
 $$\int_{\mathbb{H}^n}\f{1}{|\xi|^\beta}\zeta\le(Q,\alpha|u|^{Q^\prime}\ri)d\xi\leq \widetilde{C}$$
 for some constant $\widetilde{C}$ depending only on $C$, $Q$, and
 $r$. Hence we conclude (\ref{Ttau}).

 Secondly, we prove that for any fixed $\beta:0\leq \beta<Q$, $\alpha>0$, and
 $u\in W^{1,Q}(\mathbb{H}^n)$, there holds
 \be\label{4.13}\int_{\mathbb{H}^n}\f{1}{|\xi|_h^\beta}\zeta\le(Q,\alpha|u|^{Q^\prime}\ri)d\xi<\infty.\ee
 Since $C_0^\infty(\mathbb{H}^n)$ is dense in $W^{1,Q}(\mathbb{H}^n)$, we can take some $u_0\in C_0^\infty(\mathbb{H}^n)$
 such that $\|u-u_0\|_{W^{1,Q}(\mathbb{H}^n)}<\epsilon$, where $\epsilon>0$ is a small number to be determined later. Set
 $$w=\f{u-u_0}{\|u-u_0\|_{W^{1,Q}(\mathbb{H}^n)}}.$$
 Then $\|w\|_{W^{1,Q}(\mathbb{H}^n)}=1$.
 We divide the proof of (\ref{4.13}) into two cases:

 {\it Case $1$.} $\beta=0$.

 Recall (\ref{MF}). By (\cite{Yang}, Lemma 2.2), $\zeta(Q,t)$ is convex with respect to $t$. Since $|a+b|^{\gamma}\leq (1+\delta)|a|^{\gamma}+C(\delta,\gamma)|b|^{\gamma}$,
 $\forall a,b\in \mathbb{R},\gamma\geq 1,\delta>0$, for some constant $C(\delta,\gamma)$ depending only on $\delta$ and $\gamma$,
 we obtain
 \bna
 \int_{\mathbb{H}^n}\zeta\le(Q,\alpha|u|^{Q^\prime}\ri)d\xi&=&\int_{\mathbb{H}^n}\zeta\le(Q,\alpha|u-u_0+u_0|^{Q^\prime}\ri)d\xi\\
 &\leq&\int_{\mathbb{H}^n}\zeta\le(Q,\alpha(1+\delta)|u-u_0|^{Q^\prime}+\alpha C(\delta,Q^\prime)|u_0|^{Q^\prime}\ri)d\xi\\
 &\leq&\f{1}{\mu}\int_{\mathbb{H}^n}\zeta\le(Q,\mu\alpha(1+\delta)|u-u_0|^{Q^\prime})\ri)d\xi+
 \f{1}{\nu}\int_{\mathbb{H}^n}\zeta\le(Q,\nu\alpha C(\delta,Q^\prime)|u_0|^{Q^\prime}\ri)d\xi\\
 &\leq&\f{1}{\mu}\int_{\mathbb{H}^n}\zeta\le(Q,\mu\alpha(1+\delta)\epsilon^{Q^\prime}|w|^{Q^\prime})\ri)d\xi+
 \f{1}{\nu}\int_{\mathbb{H}^n}\zeta\le(Q,\nu\alpha C(\delta,Q^\prime)|u_0|^{Q^\prime}\ri)d\xi,
 \ena
 where ${1}/{\mu}+1/\nu=1$, $\mu>1$, $\nu>1$. Now we choose $\epsilon>0$ sufficiently small
 such that $\mu\alpha(1+\delta)\epsilon^{Q^\prime}<\alpha_Q$. By
 (\ref{Ttau}), there holds
 $$\int_{\mathbb{H}^n}\zeta\le(Q,\mu\alpha(1+\delta)\epsilon^{Q^\prime}|w|^{Q^\prime})\ri)d\xi\leq C_1$$
 for some constant $C_1$ depending only on $n$ and $\tau$. In addition, since $u_0\in C_0^\infty(\mathbb{H}^n)$, it is obvious that
 $$\int_{\mathbb{H}^n}\zeta\le(Q,\nu\alpha C(\delta,Q^\prime)|u_0|^{Q^\prime}\ri)d\xi<\infty.$$
 Therefore, we have
 $$\int_{\mathbb{H}^n}\zeta\le(Q,\alpha|u|^{Q^\prime}\ri)d\xi<\infty.$$

 {\it Case $2$.} $0<\beta<Q$.

 Note that
 \bna
 \int_{\mathbb{H}^n}\f{1}{|\xi|_h^\beta}\zeta\le(Q,\alpha|u|^{Q^\prime}\ri)d\xi\leq
 \int_{|\xi|_h\leq 1}\f{1}{|\xi|_h^\beta}\zeta\le(Q,\alpha|u|^{Q^\prime}\ri)d\xi+
 \int_{\mathbb{H}^n}\zeta\le(Q,\alpha|u|^{Q^\prime}\ri)d\xi.
 \ena
 This together with H\"older's inequality and {\it Case $1$} implies
 (\ref{4.13}).

 Finally, we confirm that for any $\alpha>\alpha_Q(1-\beta/Q)$, there holds
 $$\sup_{\|u\|_{1,\tau}\leq
1}\int_{\mathbb{H}^n}
\f{1}{|\xi|_h^\beta}\zeta\le(Q,\alpha|u|^{Q^\prime}\ri)d\xi=\infty.$$
This is based on calculations of related integrals of the Moser function sequence. We omit the details
but refer the reader to \cite{CLLY}.
 $\hfill\Box$\\

 {\bf Acknowledgements.} This work was partly supported by the NSFC
11171347 and the NCET program 2008-2011.

\end{document}